\documentclass[amscd,12pt]{amsart}
\usepackage{amscd,amsfonts,amssymb,amsmath}%,rsfs}
\numberwithin{equation}{section}
\theoremstyle{plain}
   \newtheorem{thm}{Theorem}[section]
   \newtheorem{lem}[thm]{Lemma}
   \newtheorem{prop}[thm]{Proposition}
   \newtheorem{cor}[thm]{Corollary}

   \newtheorem{rem}[thm]{Remark}
\newenvironment{pf}{\par\smallskip\noindent\emph{Proof.\;}}{\qed\par\medskip}
\newenvironment{pf*}[1]{\par\smallskip\noindent\emph{#1.\;}}{\qed\par\medskip}
\newcommand{\ab}{{\mathrm{ab}}}      % abelian
\newcommand{\adic}{{\mathrm{adic}}}  % adic
\newcommand{\cd}{{\mathrm{cd}}}      % cohomological dimension
\newcommand{\codim}{{\mathrm{codim}}}% codimension
\newcommand{\ch}{{\mathrm{ch}}}
\newcommand{\cl}{{\mathrm{cl}}}      % (cycle) class
\newcommand{\Coker}{{\mathrm{Coker}}}% cokernel
\newcommand{\et}{{\mathrm{\acute{e}t}}}    % etale
\newcommand{\Frac}{{\mathrm{Frac}}}    % Galois group
\newcommand{\Gal}{{\mathrm{Gal}}}    % Galois group
\newcommand{\h}{{\mathrm{h}}}        % henselian
\newcommand{\rH}{H}       % (co)homology group
\newcommand{\Image}{{\mathrm{Im}}}   % image
\newcommand{\loc}{{\mathrm{loc}}}    % local
\newcommand{\Nis}{{\mathrm{Nis}}}    % Nisnevich
\newcommand{\Spec}{{\mathrm{Spec}}}  % spectrum
\newcommand{\supp}{{\mathrm{Supp}}}  % support
\newcommand{\sym}{{h}}               % symbol
\newcommand{\lra}{\longrightarrow}
\newcommand{\Lra}{\Longrightarrow}
\newcommand{\hra}{\hookrightarrow}
\newcommand{\bF}{{\mathbb F}}
\newcommand{\cF}{{\mathcal F}}
\newcommand{\cI}{{\mathcal I}}
\newcommand{\cK}{{\mathcal K}}
\newcommand{\cO}{{\mathcal O}}
\newcommand{\p}{{\mathfrak p}}
\newcommand{\bZ}{{\mathbb Z}}
\newcommand{\ol}[1]{\overline{#1}}
\def\us#1#2{\underset{#1}{#2}}
\def\os#1#2{\overset{#1}{#2}}
\newcommand{\isom}{\hspace{9pt}{}^\sim\hspace{-16.5pt}\lra}
\newcommand{\leftisom}{\hspace{11pt}{}^\sim\hspace{-18.0pt}\longleftarrow}
\begin{document}

\title[]{$\ell$-adic class field theory for regular local rings
%\thanks{The research for this article was supported by JSPS Postdoctoral Fellowships for Research Abroad.}
}
%\titlerunning{$\ell$-adic class field theory for local rings}
\author[]{Kanetomo Sato}
%\institute{K. Sato \at Graduate school of Mathematics, Nagoya University\\ Furo-cho, Chikusa-ku, Nagoya 464-8602, JAPAN \\ Tel.: +81-52-789-2429\\ Fax: +81-52-789-2829\\ \email{kanetomo@math.nagoya-u.ac.jp}%\\}
\date{September 23, 2007} %Received: date / Revised: date}
\begin{abstract}
In this paper, we prove the $\ell$-adic abelian class field theory for henselian regular local rings of equi-characteristic assuming the surjectivity of Galois symbol maps, which is a $\ell$-adic variant of a result of Matsumi \cite{matsumi3}.
\end{abstract}
%\subclass{11G45 \and 19F05}
\thanks{2000 {\it Mathematics Subject Classification}$:$ 11G45 \and 19F05}
\maketitle
\section{Introduction}
Let $p$ be a prime number, and let $A$ be an excellent henselian regular local ring over $\bF_p$ of dimension $n$ with finite residue field (e.g., $\bF_p[\hspace{-1.4pt}[X_1,\dotsc,X_n]\hspace{-1.4pt}]$). Let $s$ be the closed point of $X:=\Spec(A)$ and let $D$ be an effective divisor on $X$. By standard arguments in the higher class field theory (cf.\ \cite{ks}, \cite{ras}), there is a localized Chern class map \begin{equation}\label{loc-chern} c_{X,D,\ell^r}: \rH^n_s(X_{\Nis},\cK^M_n(\cO_X,\cI_D))/\ell^r \lra \rH^{2n}_s(X_{\et},j_!\mu_{\ell^r}^{\otimes n}) \end{equation} for a prime number $\ell \ne p$ and $r>0$. Here $j$ denotes the natural open immersion $X - \supp(D) \hra X$, $\cI_D\subset \cO_X$ denotes the defining ideal of $D$ and $\cK^M_n(\cO_X,\cI_D)$ denotes a certain Milnor $K$-sheaf in the Nisnevich topology (see Notation below). We will review the construction of this map in \S\ref{sec3}. The group $\rH^n_s(X_{\Nis},\cK^M_n(\cO_X,\cI_D))$ plays the role of the $K^M_n$-id\`ele class group of $K:=\Frac(A)$ with modulus $D$. The main result of this paper is the bijectivity of this map, which is conditional with respect to the surjectivity of Galois symbol maps:
\stepcounter{thm}
\begin{thm}\label{thm:app:main}
Assume {\bf S}$(K,n,\ell)$ holds $($see Notation below for the precise statement of this condition$)$. Then $c_{X,D,\ell^r}$ is bijective.
\end{thm}
\noindent
A key to the proof of this theorem is the exactness of the cousin complex
%\begin{align*}
\[ 0  \to \rH^{n+1}(\Spec(K)_{\et},\mu_{\ell^r}^{\otimes n}) \to
  \us{y \in X^1}{\bigoplus} \ \rH^{n+2}_y(X_{\et},j_!\mu_{\ell^r}^{\otimes n})
  \to \us{y \in X^2}{\bigoplus} \ \rH^{n+3}_y(X_{\et},j_!\mu_{\ell^r}^{\otimes n})
  \to \dotsb, \]
% \to \dotsb \to \ \rH^{2n+1}_s(X_{\et},j_!\mu_{\ell^r}^{\otimes n}) \to 0, \end{align*}
which we prove using a result of Panin \cite{panin}, Theorem 5.2 (cf.\ \cite{bo}, \cite{matsumi}). In fact, we will prove this exactness for the spectrum of an arbitrary henselian local ring which is essentially \'etale over $A$ (see Lemma \ref{lem:app:key} below).
% See also Corollary \ref{lem:app:hojo} below for why we need the condition {\bf S}$(K,n,\ell)$.
Using this local-global principle, we will deduce the bijectivity of $c_{X,D,\ell^r}$ from the bijectivity of Galois symbol maps for $n$-dimensional local fields \cite{kk:cft}, p.\ 672, Lemma 14 (1) (see \S\ref{sec5} below for details).
\par
\bigskip
We explain how Theorem \ref{thm:app:main}
       is related to the class field theory of $K$.
For $D \subset X$ as before, there is a natural pairing of $\bZ/\ell^r$-modules
\stepcounter{equation}
\begin{equation}\label{pairing}
\rH^{2n}_s(X_{\et},j_!\mu_{\ell^r}^{\otimes n})
  \times \rH^1(U_{\et},\bZ/\ell^r)
    \lra \rH^{2n+1}_s(X_{\et},\mu_{\ell^r}^{\otimes n}) \isom \bZ/\ell^r.
\end{equation}
Here $U$ denotes $X - \supp(D)$, and the trace isomorphism is given by that of the closed point $s$ and the cohomological purity \cite{fujiwara}, remark after 7.1.7 (cf.\ \cite{sga4}, XIX). This pairing together with $c_{X,D,\ell^r}$ induces a map \[ \rho_{X,D,\ell^r}: \rH^n_s(X_{\Nis},\cK^M_n(\cO_X,\cI_D))/\ell^r \lra \pi_1^{\ab}(U)/\ell^r.\] Taking projective limits with respect to $D$ and $r$, we obtain the $\ell$-adic reciprocity map of the class field theory of $K$: \[ \rho_{K,\ell\text{-}\adic}: C_K:=\us{D \subset X}{\varprojlim} \ \rH^n_s(X_{\Nis},\cK^M_n(\cO_X,\cI_D)) \lra G_K^{(\ell)}, \]
where $D$ runs through all effective divisors on $X$, and $G_K^{(\ell)}$ denotes the maximal pro-$\ell$ quotient of $\Gal(K^\ab/K)$ with $K^\ab$ the maximal abelian extension of $K$.
\stepcounter{thm}
\begin{cor}\label{cor:app:main}
Assume {\bf S}$(K,n,\ell)$ holds. Then $\rho_{K,\ell\text{-}\adic}$ has dense image with respect to the Krull topology on the Galois group, and its kernel consists of elements which are divisible in $C_K$ by any power of $\ell$. \end{cor} \noindent Indeed, by a theorem of Gabber \cite{illusie}, Th\'eor\`eme 5.4 (see also \cite{riou}, Th\'eor\`eme 0.2, \cite{sga5}, I), the pairing \eqref{pairing} is a non-degenerate pairing of finite groups. This fact and Theorem \ref{thm:app:main} imply that the map $\rho_{X,D,\ell^r}$ is bijective for any $D$ and $r$. Hence $\rho_{K,\ell\text{-}\adic}$ has dense image and $\ker(\rho_{K,\ell\text{-}\adic})$ agrees with $\bigcap_{r>0}\ \ell^r \cdot C_K$.
\par Corollary \ref{cor:app:main} extends results of Saito \cite{saito} and Matsumi \cite{matsumi2} to higher-dimensional cases, and its $p$-adic variant is due to Matsumi \cite{matsumi3}.
\section*{Notation}
For a commutative ring $R$ with unity and an integer $q \ge 0$, $K^M_q(R)$ denotes the Milnor $K$-group of degree $q$ defined as $(R^{\times})^{\otimes q}/J$, where $J$ denotes the subgroup of $(R^{\times})^{\otimes q}$ generated by elements of the form $x_1 \otimes \dotsb \otimes x_q$ with $x_i + x_j=0$ or $1$ for some $1 \leq i < j \leq q$. For a local ring $R$, $R^\h$ denotes the henselization of $R$. For a field $F$ and a positive integer $m$ prime to $\ch(F)$, let $\sym^q_{F,m}$ be Galois symbol map due to Tate \cite{tate} \[ \sym^q_{F,m}:K^M_q(F)/m \lra \rH^q_{\Gal}(F,\mu_m^{\otimes q}), \] where the right hand side is the Galois cohomology of degree $q$ of the absolute Galois group of $F$. In \cite{bk}, \S5, Bloch and Kato conjectured that this map is always bijective. For a prime number $\ell\ne\ch(F)$, we define the condition {\bf S}$(F,q,\ell)$ as follows:
\begin{quote} {\bf S}$(F,q,\ell)$ : {\it $\sym^q_{F',\ell}$ is surjective for any finite separable extension $F'/F$}. \end{quote}
\begin{rem}\label{rem:app:1} $\sym^q_{F,m}$ is known to be bijective in the following cases: \begin{enumerate}
\item[$1)$] $q=0$, $q=1$ $($Hilbert's theorem $90)$, $q=2$ $($Merkur'ev-Suslin {\rm \cite{ms})}.
\item[$2)$] $F$ is a finite field. Indeed, both $K^M_q(F)$ and $\rH^q_{\Gal}(F,\mu_m^{\otimes q})$ are zero for $q \ge 2$.
\item[$3)$] $\ch(F) \not = 2$ and $m=2^r$ $($Merkur'ev-Suslin {\rm \cite{ms2}}, Voevodsky {\rm \cite{voe})}.
\item[$4)$] $F$ is a henselian discrete valuation ring of characteristic zero with residue field of characteristic $p > 0$ and $m=p^r$ $($Bloch-Kato {\rm \cite{bk})}.
\item[$5)$] The case that $F$ is a henselian discrete valuation ring with residue field $k$, that $m$ is prime to $\ch(k)$ and that $\sym^q_{k,m}$ and $\sym^{q-1}_{k,m}$ are bijective $(${\rm \cite{kk:cft}}, {\rm p.\ 672}, Lemma $14$ $(1))$. In particular, in view of $(2)$ and $(4)$, $\sym^q_{F,m}$ is bijective for any local field $F$ in the sense of loc.\ cit., and any $m$ prime to $\ch(F)$. The case $\ch(F)>0$ of this fact will be useful later.\end{enumerate} \end{rem}
\par For a scheme $X$, $\cK^M_q(\cO_X)$ denotes the Milnor $K$-sheaf of degree $q$ in the Nisnevich topology on $X$. For an ideal sheaf $\cI \subset \cO_X$, we define \[ \cK^M_q(\cO_X,\cI) := \ker\big(\cK^M_q(\cO_X) \to \cK^M_q(\cO_X/\cI)\big) \] For $q \ge 0$, $X^q$ denotes the set of all points on $X$ of codimension $q$. For a point $x \in X$, $\kappa(x)$ denotes the residue field of the local ring $\cO_{X,x}$.
\section{Preliminary on henselian discrete valuation fields}
\label{sec2}
Let $F$ be a henselian discrete valuation field and let $O_F$ be its integer ring. Let $\pi$ be a prime element of $O_F$. Put $X:=\Spec(O_F)$ and denote the generic point (resp.\ closed point) of $X$ by $\eta$ (resp.\ $x$). Let $U^1_F=\{1+\pi\, a \; |\ a \in O_F\}$ be the first unit group of $O_F$, and let $U^1K^M_q(F)$ be the subgroup of $K^M_q(F)$ generated by the image of $U^1_F \otimes (F^{\times})^{\otimes q-1}$. The following elementary facts will be useful later:
\begin{lem}\label{lem:app:2} Let $m$ be a positive integer invertible in $O_F$. \begin{enumerate}
% \item[(1)] The first unit group $U^1_L=\{1+\pi\, a \; |\ a \in O_L\} \subset O_L^{\times}$ is contained in $(O_L^{\times})^m$.
\item[(1)] Let $I$ be a non-zero ideal of $O_F$ with $I \ne O_F$, and let $K^M_q(O_F,I)$ be the kernel of the natural map $K^M_q(O_F) \to K^M_q(O_F/I)$. Then we have \[ \Image\big(K^M_q(O_F,I) \to K^M_q(F)\big) \subset U^1K^M_q(F) \subset m \cdot K^M_q(F). \]
\item[(2)] The connecting map $\rH^q(\eta_{\et},\mu_m^{\otimes q}) \to \rH^{q+1}_x(X_{\et},\mu_m^{\otimes q})$ in the localization theory of \'etale cohomology is surjective for any $q \ge 1$. \end{enumerate} \end{lem}
\begin{pf}
(1) The second inclusion follows from the fact that $U^1_F \subset (O_F^{\times})^m$, which is a direct consequence of Hensel's lemma. We show the first inclusion. Let $k$ be the residue field of $O_F$. There is a short exact sequence (cf.\ \cite{kk:cft}, p.\ 616, Lemma 6) \[ 0 \lra U^1K^M_q(F) \lra K^M_q(F) \os{\partial}{\lra} K^M_q(k) \oplus K^M_{q-1}(k) \lra 0, \] where the arrow $\partial$ is defined by the assignment \[\begin{CD} \{\pi,a_1,\dotsc,a_{q-1} \} @.~ \mapsto ~@. \big( @. 0 @. , @. \;\{\ol{a_1},\dotsc,\ol{a_{q-1}}\} @. \big) \\ \{a_1,\dotsc,a_{q} \} @. \mapsto @. \big(@. \{\ol{a_1},\dotsc,\ol{a_{q}}\} @. , @. 0 @. \big) \end{CD}\] with each $a_i \in O_F^\times$ and $\ol {a_i}$ denotes the residue class of $a_i$ in $k^\times$. The assertion follows from this exact sequence and the assumption that $I \ne O_F$.
% (2) follows from (1) and the fact that $\ker(\partial)\subset K^M_m(F)$ is generated by the image of $K^M_m(O_F)$ and symbols of the form $\{a,b_1,\dotsc,b_{m-1}\}$ with $a \in U^1_F$ and each $b_i \in F^\times$ (cf.\ \cite{BT}, I.4.3).
\par (2) There is a natural map \[ \rH^{q-1}(x_{\et},\mu_m^{\otimes q-1}) \lra \rH^{q+1}_x(X_{\et},\mu_m^{\otimes q}) \] sending $\alpha \in \rH^{q-1}(x_{\et},\mu_m^{\otimes q-1})$ to $\cl_X(x) \cup \alpha$, where $\cl_X(x) \in \rH^2_{x}(X_{\et},\mu_m)$ denotes the cycle class of $x$ (\cite{sga4.5}, Cycle, \S2.1). This map is bijective by the purity for discrete valuation rings (\cite{sga5}, I.5). Let $\alpha \in \rH^{q-1}(x_{\et},\mu_m^{\otimes q-1})$ be an arbitrary cohomology class. Since $X$ is henselian local with closed point $x$, there is a unique $\alpha' \in \rH^{q-1}(X_{\et},\mu_m^{\otimes q-1})$ that lifts $\alpha$. By \cite{sga4.5}, Cycle, 2.1.3, the cohomology class $\{ \pi \} \cup \alpha' \in \rH^q(\eta_{\et},\mu_m^{\otimes q})$ maps to $\alpha$ under the composite map \[\rH^q(\eta_{\et},\mu_m^{\otimes q}) \lra \rH^{q+1}_x(X_{\et},\mu_m^{\otimes q}) \leftisom \rH^{q-1}(x_{\et},\mu_m^{\otimes q-1}),\]up to a sign, which implies the assertion.
\end{pf}
\section{Localized Chern class maps}\label{sec3}
In this section, we construct localized Chern class maps.
Our construction is essentially the same as those in \cite{ks}, \cite{ras}.
We work under the following setting. Let $X$ be a noetherian integral normal scheme. Let $D$ be an effective Weil divisor on $X$. Let $j$ be the open immersion $X - \supp(D) \hra X$ and let $\cI_D\subset \cO_X$ be the defining ideal of $D$. Let $\ell$ be a prime number invertible on $X$, and let $r$ be a positive integer. For $q \ge 0$ and $x \in X^a$, we define the localized Chern class map
\begin{equation}\label{loc-chern2}
\cl_{X,D,x,\ell^r}^{q,\loc}:\rH^a_x(X_{\Nis},\cK_q^{M}(\cO_X,\cI_D))/\ell^r \lra \rH^{q+a}_x(X_{\et},j_!\mu_{\ell^r}^{\otimes q})
\end{equation}
by induction on $a \geq 0$. The map $c_{X,D,\ell^r}$ in \eqref{loc-chern} is defined as $\cl_{X,D,s,\ell^r}^{n,\loc}$. We first recall the following general facts.
\stepcounter{thm}
\begin{lem}\label{lem:app:1}
Let $Z$ be a scheme of dimension $d$, and let $\cF$ be an abelian sheaf on $Z_{\Nis}$. Then$:$
\begin{enumerate}
\item[$1)$]
We have $\rH^t(Z_{\Nis},\cF)=0$ for any $t>d$.
\item[$2)$]
For $z \in Z$, $\rH^t_z(Z_{\Nis},\cF)$ is isomorphic to
\[\begin{cases}
 \Coker\big(\cF_z \to \rH^0(U_{\Nis},\cF)\big) & \quad \hbox{ if } t=1 \\
 \rH^{t-1}(U_{\Nis},\cF) & \quad \hbox{ if } t \geq 2 \\
 0 & \quad \hbox{ if } t \geq \codim_Z(z)+1, \\
\end{cases}\]
where $U$ denotes $\Spec(\cO_{Z,z}^{\h})-\{ z \}$.
\item[$3)$]
There is an exact sequence
$$
\us{z \in Z^{d-1}}{\bigoplus} \rH^{d-1}_z(Z_{\Nis},\cF)
   \lra
   \us{z \in Z^d}{\bigoplus} \rH^d_z(Z_{\Nis},\cF)
      \lra \rH^d(Z_{\Nis},\cF)
        \lra 0.
$$
\end{enumerate}
\end{lem}
\begin{pf}
For (1), see \cite{ras}, Lemma 1.22. The assertion (2) follows from (1) and a standard localization argument using the excision in the Nisnevich topology. To show (3), we compute the local-global spectral sequence \[ E_1^{u,v}=\us{z \in Z^u}{\bigoplus} \rH^{u+v}_z(Z_{\Nis},\cF) \Lra \rH^{u+v}(Z_{\Nis},\cF). \] We have $E_1^{u,v}=0$ unless $0 \le u \le d$, and $E_1^{u,v}=0$ for $v>0$ by (2). Hence we have $E_2^{d,0}\simeq \rH^d(Z_{\Nis},\cF)$, which implies (3).
\end{pf}
We construct the map \eqref{loc-chern2} in three steps. Replacing $X$ by $X_x^{\h}$, the henselization of $X$ at $x$ if necessary, we assume that $X$ is henselian local with closed point $x$. % Put $a:=\codim_X(x)(=\dim(X))$.
\par
\vspace{8pt}
\noindent
{\bf Step 0.} \hspace{1pt} Assume $a=0$. Then we have $X=x$ and $D$ is zero. We define $\cl_{X,0,X,\ell^r}^{q,\loc}$ as the Galois symbol map $\sym_{x,\ell^r}^q: K^M_q(x)/\ell^r \to \rH^q(x_{\et},\mu_{\ell^r}^{\otimes q})$.
\par
\vspace{8pt}
\noindent
{\bf Step 1.} \hspace{1pt} Assume $a=1$. Then $X$ is the spectrum of a henselian discrete valuation ring $O_F$ with fraction field $F$. Let $\eta=\Spec(F)$ be the generic point of $X$. We have $\supp(D)= x$, or otherwise, $D$ is zero. Put $\cK_D:= \cK_q^{M}(\cO_X,\cI_D)$ and $\mu_D:= j_!\mu_{\ell^r}^{\otimes q}$ for simplicity, which mean $\cK_q^{M}(\cO_X)$ and $\mu_{\ell^r}^{\otimes q}$, respectively, if $D=0$ (i.e., $\cI_D=\cO_X$). We define $\cl_{X,D,x,\ell^r}^{m,\loc}$ as the map induced by the commutative diagram with exact rows
\stepcounter{equation}
\begin{equation}\label{cd:app:0}
\begin{CD}
\rH^0(X_{\Nis},\cK_D)/\ell^r @. \; \lra \; @. K^M_q(F)/\ell^r @. \; \lra \; @. \rH^1_x(X_{\Nis},\cK_D)/\ell^r @. \; \lra \; @. 0 \phantom{.} \\ @VVV @. @V{\sym_{\eta,\ell^r}^q}VV \\ \rH^q(X_{\et},\mu_D) @. \; \lra \; @. \rH^q(\eta_{\et},\mu_{\ell^r}^{\otimes q}) @. \; \os{\delta}{\lra} \; @. \rH^{q+1}_x(X_{\et},\mu_D) @. \; \lra \; @. 0. \\
\end{CD}
\end{equation}
Here the upper row (resp.\ lower row) arises from a localization long exact sequence of Nisnevich (resp.\ \'etale) cohomology groups, and the upper row is exact by Lemma \ref{lem:app:1} (1). The arrow $\delta$ is surjective by Lemma \ref{lem:app:2} (2) (resp.\ the fact that $\mu_D|_x=0$) if $D=0$ (resp.\ $D \ne 0$). The left vertical arrow is the Galois symbol map $K^M_q(O_F)/\ell^r \to \rH^q(X_{\et},\mu_{\ell^r}^{\otimes q})$ (resp.\ zero map) if $D=0$ (resp.\ $D \ne 0$). The left square commutes obviously if $D=0$. If $D \ne 0$, then the top left arrow is zero by Lemma \ref{lem:app:2} (1), and the left square commutes. Thus we obtain the map $\cl_{X,D,x,\ell^r}^{q,\loc}$ in the case $a=1$.
\par
\vspace{8pt}
\noindent
{\bf Step 2.}
\hspace{1pt}
Assume $a \ge 2$ and that the localized Chern class maps have been defined for points of codimension $\le a-1$. Put $\cK_D:= \cK_q^{M}(\cO_X,\cI_D)$ and $\mu_D:= j_!\mu_{\ell^r}^{\otimes q}$ as before. We define $\cl_{X,D,x,\ell^r}^{q,\loc}$ to be the map induced by the commutative diagram
\begin{equation}{\small \begin{CD} \us{y \in X^{a-2}}{\bigoplus} \rH^{a-2}_y(X_{\Nis},\cK_D)/\ell^r @. \; \to \; @. \us{y \in X^{a-1}}{\bigoplus} \rH^{a-1}_y(X_{\Nis},\cK_D)/\ell^r @. \; \to \; @. \rH^a_x(X_{\Nis},\cK_D)/\ell^r @. \; \to \; @. 0\\ @V{\bigoplus \cl_{X,D,y,\ell^r}^{q,\loc}}VV @. @V{\bigoplus \cl_{X,D,y,\ell^r}^{q,\loc}}VV @. @. @. \\ \us{y \in X^{a-2}}{\bigoplus} \rH^{q+a-2}_y(X_{\et},\mu_D) @.\; \to \; @. \us{y \in X^{a-1}}{\bigoplus} \rH^{q+a-1}_y(X_{\et},\mu_D) @. \; \to \; @. \rH^{q+a}_x(X_{\et},\mu_D),\end{CD}}
\label{cd:app:4}\end{equation}
where the upper row is exact by Lemma \ref{lem:app:1} (3) for $X - \{x \}$ and the lower row is a complex. The left square commutes by the construction of the vertical arrows (cf.\ \eqref{cd:app:0}). This completes the construction of the localized Chern class map \eqref{loc-chern2}.
\addtocounter{thm}{2}
\begin{rem}
One can check that the group $\rH^a_x(X_{\Nis},\cK^M_q(\cO_X,\cI_D))/\ell^r$ depends only on $X - \supp(D)$ by repeating the arguments in this section.
\end{rem}
\section{Key results on \'etale cohomology}
\label{sec4}
Let $A$ be as in the introduction and let $\ell$ be a prime number invertible in $A$. In this section we prove Lemma \ref{lem:app:key} below (compare with Lemma \ref{lem:app:1}) using the following key fact due to Gabber (\cite{illusie}, Corollaire 4.3, \cite{orgogozo}, Th\'eor\`eme 5.1, \cite{PS}, Th\'eor\`eme 4):
\begin{thm}[{\bf Gabber}]\label{thm:gabber}
Let $R'$ be an integral excellent henselian local ring of dimension $d$ with residue field $k$, and let $K'$ be the fraction field of $R'$. Then we have $\cd_\ell(K') = d + \cd_\ell (k)$ for any prime number $\ell\ne\ch(k)$.
\end{thm}
\begin{lem}\label{lem:app:key}
Let $R$ be a henselian local ring which is essentially \'etale over $A$ $(R$ may be $A$ itself\,$)$. Put $X:=\Spec(R)$ and $a:=\dim(X)$. Then$:$
\begin{enumerate}
\item[$1)$]
For $y \in X$, we have $\cd_{\ell}(y) \leq n+1-\codim_X(y)$.
\item[$2)$]
For an $\ell$-primary torsion sheaf $\cF$ on $X_{\et}$ and $y \in X^c$ $(c \ge 0)$, $\rH^t_y(X_{\et},\cF)$ is isomorphic to
\[ \begin{cases}
    \Coker\big(\rH^{t-1}(\Spec(\cO_{X,y}^{\h})_{\et},\cF)
       \to \rH^{t-1}(U_{\et},\cF) \big)
        & \hbox{ if \;} t= n+2-c \\
    \rH^{t-1}(U_{\et},\cF)
        & \hbox{ if \;} t \geq n+3-c \\
    0   & \hbox{ if \;} t \geq n+2+c, \\
   \end{cases} \] where $U$ denotes $\Spec(\cO_{X,y}^{\h})-\{ y \}$.
\item[$3)$]
Assume $a \ge 1$. Let $x$ $($resp.\ $\eta)$ be the closed point $($the generic point$)$ of $X$. Let $D$ be an effective divisor on $X$. Let $j$ be the open immersion $X - \supp(D) \hra X$ and put $\mu_D:=j_!\mu_{\ell^r}^{\otimes n}$. Then the sequence
\begin{align*}
0  \to \rH^{n+1}(\eta{\,}_{\et},\mu_{\ell^r}^{\otimes n}) \to \us{y \in X^1}{\bigoplus} \ \rH^{n+2}_y(X_{\et},\mu_D) \to \ & \us{y \in X^2}{\bigoplus} \ \rH^{n+3}_y(X_{\et},\mu_D) \to \\ \dotsb \to \ & \rH^{n+a+1}_{x}(X_{\et},\mu_D) \to 0
\end{align*}
is exact.
\end{enumerate}
\end{lem}
\noindent
Admitting this lemma,
we first prove the following two consequences, where the notation remains as in Lemma \ref{lem:app:key} and $K$ denotes $\Frac(A)$.
\begin{cor}\label{cor:app:key}
If $a \ge 2$, then the sequence
$$
    \us{y \in X^{a-2}}{\bigoplus}\rH^{n+a-2}_y(X_{\et},\mu_D)
    \to
    \us{y \in X^{a-1}}{\bigoplus}\rH^{n+a-1}_y(X_{\et},\mu_D)
    \to
    \rH^{n+a}_{x}(X_{\et},\mu_D)
    \to 0
$$
is exact.
\end{cor}
\begin{pf}
Put $Y := X - \{x \}$. Since $n+3-a<n+a$, we have $\rH^{n+a-1}(Y_{\et},\mu_D) \simeq \rH^{n+a}_{x}(X_{\et},\mu_{D})$ by Lemma \ref{lem:app:key} (2).
We prove the sequence
\addtocounter{equation}{3}
\begin{equation}\label{seq:app:1}
\us{y \in Y^{a-2}}{\bigoplus} \rH^{n+a-2}_y(Y_{\et},\mu_D)
   \to \us{y \in Y^{a-1}}{\bigoplus} \rH^{n+a-1}_y(Y_{\et},\mu_D)
      \to \rH^{n+a-1}(Y_{\et},\mu_D)
        \to 0
\end{equation}
is exact by computing the spectral sequence
\begin{equation}\label{ss:loc:Y}
 E_1^{u,v}=\bigoplus_{y \in Y^u}
    \rH^{u+v}_y(Y_{\et},\mu_D)
       \Lra \rH^{u+v}(Y_{\et},\mu_D).
\end{equation}
Since $\dim(Y)=a-1$, $E_1^{u,v}$ is zero unless $0 \le u \le a-1$. We have $E_1^{u,v}=0$ for $v \geq n+2$ by Lemma \ref{lem:app:key} (2), and $E_2^{u,n+1}$ is zero for $0 \leq u \leq a-2$ by Lemma \ref{lem:app:key} (3) for $X$. Hence we obtain $E_2^{a-1,n} \simeq \rH^{n+a-1}(Y_{\et},\mu_D)$ and the exactness in question.
\end{pf}
\addtocounter{thm}{2}
\begin{cor}\label{lem:app:hojo}
Assume {\bf S}$(K,n,\ell)$ holds. Then the map \[\cl_{X,D,x,\ell^r}^{n,\loc}: \rH^{a}_x(X_{\Nis},\cK_n^{M}(\cO_X,\cI_D))/\ell^r \lra \rH^{n+a}_x(X_{\et},j_!\mu_{\ell^r}^{\otimes n})\] is surjective.
\end{cor}
\begin{pf}
If $a \le 1$, then $\cl_{X,D,X,\ell^r}^{n,\loc}$ is surjective by its construction in \S\ref{sec3} and the assumption {\bf S}$(K,n,\ell)$ (see also \cite{bk}, 5.13 (i)). Here we have used the surjectivity of the lower $\delta$ in \eqref{cd:app:0}. We prove the case $a \geq 2$ by induction on $a$. By the induction hypothesis and the construction of $\cl_{X,D,x,\ell^r}^{n,\loc}$ (cf.\ \eqref{cd:app:4}), it is enough to show that the connecting map \[ \us{y \in X^{a-1}}{\bigoplus} \rH^{n+a-1}_y(X_{\et},\mu_{D}) \lra \rH^{n+a}_x(X_{\et},\mu_{D}) \] is surjective, which has been shown in Corollary \ref{cor:app:key}.
\end{pf}
In the rest of this section, we prove Lemma \ref{lem:app:key}.
\par\smallskip
\begin{pf*}{Proof of Lemma $\mathrm{\ref{lem:app:key}}$}
(1) By \cite{se}, I.3.3, Proposition 14, it is enough to deal with the case $R=A$. Then the assertion follows from Theorem \ref{thm:gabber} applied to $R'=\Image(A \to \kappa(y))$ and the assumption that the residue field of $A$ is finite. \par (2) follows from (1) and a standard localization argument using the excision in the \'etale topology. The details are straight-forward and left to the reader.
%
% By \cite{se}, I.3.3, Proposition 14, it suffices to show the case $R=A$. We prove this case by induction on $n$. The case $n=0$ (i.e., $A = \bF_q$) is obvious. Suppose that Lemma \ref{lem:app:key} (1) has been proved for any formal power series ring over $\bF_q$ of dimension $\leq n-1$. Let $y$ be a point on $X=\Spec(A)$ $(A=\bF_q[\hspace{-1.4pt}[T_1,\cdots,T_n]\hspace{-1.4pt}])$, and put $c:=\codim_X(x)$. Let $\p_y$ be the prime ideal of $A$ corresponding to $y$. If $c=0$, then $y$ is the spectrum of the fraction field of $A$, and we have $\cd_{\ell}(y) \leq n+1$ by a theorem of Gabber included in \cite{matsumi}. Next suppose that $c \geq 1$. We take a prime ideal $\p' \subset A$ of height $1$ with $\p' \subset \p_y$. Then by Weierstrass' preparation theorem, the domain $A/\p'$ is a finite integral extension over a formal power series ring $\bF_q[\hspace{-1.4pt}[t_1,\cdots,t_{n-1}]\hspace{-1.4pt}]$. The prime ideal $\p_y/\p' \subset A/\p'$ has height $c-1$, and the prime ideal ${\mathfrak q} := (\p_y/\p') \cap \bF_q[\hspace{-1.4pt}[t_1,\cdots,t_{n-1}]\hspace{-1.4pt}]$ of $\bF_q[\hspace{-1.4pt}[t_1,\cdots,t_{n-1}]\hspace{-1.4pt}]$ also has height $c-1$ by the going-up and going-down theorems \cite{am}, 5.11, 5.16. Hence we have $\cd_{\ell}(y) \leq \cd_{\ell}(\kappa({\mathfrak q})) \leq n+1-c$ by \cite{se}, I.3.3, Proposition 14 and the induction hypothesis.
%
\par (3) By a result of Panin \cite{panin}, 5.2 (cf.\ \cite{bo}), the sequence
%      with $\mu_{\ell^r}^{\otimes n}$-coefficients
\addtocounter{equation}{1}
\begin{equation}\label{app:cousin:ex}
\begin{CD}
0 \to \rH^{n+1}(\eta{\,}_{\et},\mu_{\ell^r}^{\otimes n}) \to \us{y \in X^1}{\bigoplus}\rH^{n+2}_y(X_{\et},\mu_{\ell^r}^{\otimes n}) \to \us{y \in X^2}{\bigoplus}\rH^{n+3}_y(X_{\et},\mu_{\ell^r}^{\otimes n}) \to \dotsb \\ \to \rH^{n+a+1}_x(X_{\et},\mu_{\ell^r}^{\otimes n}) \to 0
\end{CD}\end{equation} is exact, where we have used the vanishing $\rH^{n+1}(X_{\et},\mu_{\ell^r}^{\otimes n})=0$. Indeed, we have \[ \rH^q(X_{\et},\mu_{\ell^r}^{\otimes n}) \simeq \rH^{q}(x_{\et},\mu_{\ell^r}^{\otimes n}) =0 \qquad \hbox{ for } q>n+1-a \] by Lemma \ref{lem:app:key} (1), and we have $n+1>n+1-a$ by the assumption that $a(=\dim(X)) \geq 1$. Thus we obtain the assertion if $D=0$. In particular,
\begin{equation}\label{app:isom:codim1}
  \rH^{n+1}(\eta{\,}_{\et},\mu_{\ell^r}^{\otimes n})
     \simeq \rH^{n+2}_x(X_{\et},\mu_{\ell^r}^{\otimes n}),
    \qquad \hbox{ if } a=1.
\end{equation}
In what follows, assume $D \ne 0$. We prove the sequence in question is isomorphic to the exact sequence \eqref{app:cousin:ex}. It is enough to show the natural map
\begin{equation}\label{app:isom:codima}
  \rH^{n+a+1}_x(X_{\et},\mu_D) \lra \rH^{n+a+1}_x(X_{\et},\mu_{\ell^r}^{\otimes n})
\end{equation}
is bijective for any $X$ and any non-zero $D$ on $X$ as in Lemma \ref{lem:app:key}, where $x$ denotes the closed point of $X$ and $a$ denotes $\dim(X)$. We prove this assertion by induction on $a=\dim(X)$. If $a=1$, then the assertion follows from the isomorphisms \[ \rH^{n+2}_x(X_{\et},\mu_D) \leftisom \rH^{n+1}(\eta{\,}_{\et},\mu_{\ell^r}^{\otimes n}) \us{\eqref{app:isom:codim1}}{\isom} \rH^{n+2}_x(X_{\et},\mu_{\ell^r}^{\otimes n}), \] where the left arrow is bijective by the fact that $\mu_D\vert_x=0$. Assume $a \geq 2$. Computing the spectral sequence \eqref{ss:loc:Y} using Lemma \ref{lem:app:key} (2), one can easily check the sequence \[ \us{y \in Y^{a-2}}{\bigoplus} \rH^{n+a-1}_y(Y_{\et},\mu_D) \to \us{y \in Y^{a-1}}{\bigoplus} \rH^{n+a}_y(Y_{\et},\mu_D) \to \rH^{n+a}(Y_{\et},\mu_D) \to 0 \] is exact. Since $\mu_D\vert_x=0$, we have $\rH^{n+a}(Y_{\et},\mu_D) \simeq \rH^{n+a+1}_x(X_{\et},\mu_D)$. Thus the sequence \[ \us{y \in X^{a-2}}{\bigoplus} \rH^{n+a-1}_y(X_{\et},\mu_D) \to \us{y \in X^{a-1}}{\bigoplus} \rH^{n+a}_y(X_{\et},\mu_D) \to \rH^{n+a+1}_x(X_{\et},\mu_D) \to 0 \] is exact. We obtain the bijectivity of \eqref{app:isom:codima} by comparing this exact sequence with \eqref{app:cousin:ex} using the induction hypothesis. This completes the proof of the lemma.
\end{pf*}
\section{Proof of Theorem \ref{thm:app:main}}
\label{sec5}
We introduce some auxiliary terminology and notation.
For $0 \leq a \leq n-1$, we define a chain over $A$ of length $a$ to be a sequence
\begin{equation}\label{app:chain}
\big(
{\mathfrak m}_A, \p_{n-1}^n,\p_{n-2}^{n-1},\p_{n-3}^{n-2},\cdots,\p_{n-a}^{n-a+1}
\big),
\end{equation}
where ${\mathfrak m}_A$ is the maximal ideal of $A$ and $\p_{n-1}^n$ is a prime ideal of $A$ of height $n-1$. For $2 \leq q \leq a$, $\p_{n-q}^{n-q+1}$ is a prime ideal of height $n-q$ of the $(n-q+1)$-dimensional henselian local ring
$$
\left( \cdots \left( \left(
   A^{\h}_{\p_{n-1}^n} \right){}^{\h}_{\p_{n-2}^{n-1}} \right)
       \cdots \right){}^{\h}_{\p_{n-q+1}^{n-q+2}}.
$$
For a chain $\delta$ over $A$ of length $a$ of the form \eqref{app:chain},
   we define the henselian local ring $A_{\delta}^{\h}$
     of dimension $n-a$ as
$$
A_{\delta}^{\h} := \left( \cdots \left( \left(
   A^{\h}_{\p_{n-1}^n} \right){}^{\h}_{\p_{n-2}^{n-1}} \right)
       \cdots \right){}^{\h}_{\p_{n-a}^{n-a+1}}.
$$
For $\delta$ of length zero, $A_{\delta}^{\h}$ means $A$ itself.
We prove
\addtocounter{thm}{1}
\begin{prop}\label{prop:app:key}
Assume $1 \leq a \leq n$ and that {\bf S}$(K,n,\ell)$ holds. Let $\delta$ be a chain over $A$ of length $n-a$, and let $D$ be an effective divisor on $X:= \Spec(A_{\delta}^{\h})$. Let $j$ be the open immersion $X - \supp(D) \hra X$ and let $x$ be the closed point of $X$. Then the map
\[ \cl_{X,D,x,\ell^r}^{n,\loc}: \rH^{a}_x(X_{\Nis},\cK_n^{M}(\cO_X,\cI_D))/\ell^r \lra \rH^{n+a}_x(X_{\et},j_!\mu_{\ell^r}^{\otimes n}) \] is bijective.
\end{prop}
\noindent
Theorem \ref{thm:app:main}
 follows from the case $a=n$ of Proposition \ref{prop:app:key}.
\par \smallskip
\begin{pf}
Assume $a=1$. The generic point $\eta$ (resp.\ the closed point $x$) of $X$ is the spectrum of an $n$-dimensional local field (resp.\ $(n-1)$-dimensional local field) in the sense of \cite{kk:cft}, and the Galois symbol maps of $\eta$ and $x$ are bijective by Remark \ref{rem:app:1} (5). Hence $\cl_{X,D,x,\ell^r}^{n,\loc}$ is bijective by the diagram \eqref{cd:app:0} and Lemma \ref{lem:app:2} (1).
%, noting that the top left arrow in \eqref{cd:app:1} is zero by \eqref{incl:app:1}. We proceed our proof by induction on $a \ge 1$.
The case $a \ge 2$ follows from the induction on $a$ using the diagram \eqref{cd:app:4} and Corollaries \ref{cor:app:key} and \ref{lem:app:hojo}.
\end{pf}
\noindent
This completes the proof of Theorem \ref{thm:app:main}.
\par
\bigskip
\noindent
%\begin{acknowledgements}
{\bf Acknowledgements.} The author expresses his gratitude to Professor Luc Illusie for valuable comments on the references for Gabber's results (\cite{illusie}, \cite{orgogozo}, \cite{PS}, \cite{riou}), and to Atsushi Shiho and Kazuya (Pierre) Matsumi for stimulating discussions.
 % Professor Wayne Raskind,
%\end{acknowledgements}

\end{document}